\documentclass [12pt]{article}
\usepackage{amsmath,amssymb}

\begin{document}

   \title {The column and row immanants of matrices
 over a split quaternion  algebra.}

\author {Ivan Kyrchei \footnote{Pidstrygach Institute
for Applied Problems of Mechanics and Mathematics, str.Naukova 3b,
Lviv, Ukraine, 79060, kyrchei@lms.lviv.ua}}
\date{}
\maketitle
\begin{abstract}
The theory of the column-row determinants has been considered for matrices over a non-split quaternion algebra.
In this paper the concepts of  column-row determinants are  extending  to a split quaternion algebra. New definitions of the column and row immanants (permanents) for matrices over a non-split quaternion algebra are introduced, and their basic  properties are investigated.
  The key theorem about the column and row immanants of a Hermitian matrix over a split quaternion algebra is proved. Based on this theorem an  immanant of a Hermitian matrix over a split quaternion algebra is introduced.
\end{abstract}

\section{Introduction}

\newtheorem{thm}{Theorem}[section]
\newtheorem{lem}{Lemma}[section]
\newtheorem{defn}{Definition}[section]

\newcommand{\rank}{\mathop{\rm rank}\nolimits}
\newtheorem{prop}{Proposition}[section]
The immanent of a matrix  is a generalization of the concepts of determinant and permanent. The immanent of a complex matrix  was defined by Dudley E. Littlewood and Archibald Read Richardson in \cite{li} as follows.
 \begin{defn}
 Let  $\sigma \in S_{n}$ denote the symmetric group on $n$ elements. Let $\chi: S_{n}\rightarrow {\mathbb{C}}$  be a complex character. For any $n\times n$ matrix ${\rm {\bf A}}=(a_{ij})\in{\mathbb{C}}^{n \times n} $ define the immanent of ${\rm {\bf A}}$ as
\[
{\rm{Imm}}_{\chi}( {\bf A})={\sum\limits_{\sigma \in S_{n}}}\chi(\sigma)\prod_{i=1}^{n}a_{i\,\sigma(i)}
\]
\end{defn}
Special cases of immanants are determinants and permanents. In the case where $\chi$ is the constant character $(\chi(x)=1$ for all $x\in S_{n}$), ${\rm{Imm}}_{\chi}( {\bf A})$ is the permanent of ${\rm {\bf A}}$. In the case where $\chi$ is the ${\rm sign}$ of the permutation (which is the character of the permutation group associated to the (non-trivial) one-dimensional representation), ${\rm{Imm}}_{\chi}( {\bf A})$ is the determinant of ${\rm {\bf A}}$.
The main goal of this paper is the extending  the concept of  immanant  to a split quaternion algebra using methods of the theory of the row and column determinants. The theory of the row and column determinants was introduced
 in \cite{ky1, ky2} for matrices   over the quaternion non-split algebra.
This theory  over the quaternion skew field is being actively developed as by  the author \cite{ky3}-\cite{ ky5}, and others (see, for ex. \cite{so1}-\cite{so3}).

The paper is organized as follows. In Section
2 we consider briefly the main provisions of the quaternion algebra.
In Section 3 definitions of the row and column immanents (consequently, determinants and permanents) are given. Their properties of an arbitrary quadratic matrix over the quaternion algebra
 are described in Section 4.
In Section 5 the key theorem about the column and row immanants of a Hermitian matrix over a split quaternion algebra is proved and based on it we introduce the immanant (determinant, permanent) of a Hermitian matrix.

\section{Quaternion algebra}
 A quaternion algebra ${\bf{H}}( a,b)$ over a field ${\bf{F}}$ is a central simple algebra over ${\bf{F}}$ that is a four-dimensional vector space over
${\bf{F}}$  with basis $\{1,i,j,k\}$
and the following multiplication rules:
\[
  i^{2}=a, \,\,\,
  j^{2}=b, \,\,\,
  ij=k, \,\,\,
  ji=-k.
\]
A quaternion algebra ${\bf{H}}( a,b)$ over ${\bf{F}}$ is denoted $(\frac{\alpha,\beta}{{\bf{F}}})$ as well. To every quaternion algebra ${\bf{H}}( a,b)$, one can associate a quadratic form ${\rm n}$
(called the norm form) on ${\bf{H}}$ such that ${\rm
n}(xy)={\rm n}(x){\rm n}(y)$ for all $x$ and $y$ in
${\bf{H}}$. A linear mapping $x \rightarrow \overline{x}={\rm
t}(x)-x$ is also defined on ${\bf{H}}$. It is an involution,
i.e. $\overline{\overline{x}}=x$, $\overline {x + y} = \overline
{x} + \overline {y} $ and $\overline {x \cdot y} = \overline {y}
\cdot \overline {x} $
 An element
$\overline{x}$ is called the conjugate of $x \in{\bf{H}}$. ${\rm t}(x)$ and ${\rm n}(x)$ are called the trace
and the norm of $x$ respectively, at that $\{{\rm n}(x), {\rm
t}(x)\} \subset {\bf{F}}$ for all $x$  in ${\bf{H}}$. They
also satisfy the following conditions:  ${\rm n}\left(
\overline{x} \right) = {\rm n}(x)$, ${\rm t}\left( \overline{x}
\right) = {\rm t}(x)$ and ${\rm t}\left( {q \cdot p} \right) =
{\rm t}\left( {p \cdot q} \right)$. The last property is the
rearrangement property of the trace.

Depending on the choice of ${\bf{F}}$, $ a$ and $b$ we
have only two possibilities (\cite{le}):

1.$(\frac{a,b}{{\bf{F}}})$ is a division algebra,

2. $(\frac{a,b}{{\bf{F}}})$ is isomorphic to
 the algebra of all $2\times 2$ matrices with entries from ${\bf{F}}$.

 If an ${\bf{F}}$-algebra is isomorphic to a full matrix algebra over ${\bf{F}}$ we say
that the algebra is split, so (2) is the split case.

The most famous example of a non-split quaternion algebra is Hamilton's quaternions  ${\rm {\mathbb{H}}}=(\frac{-1,-1}{{\rm {\mathbb{R}}}})$.

An example of a split quaternion algebra is  split quaternions of James Cockle ${\bf{H}}_{\bf{S}}(\frac{-1,1}{{\rm {\mathbb{R}}}})$.  Recently there was conducted a number of studies in split quaternion matrices (see, for ex. \cite{er1}-\cite{al}). The matrix representation for the complex
quaternions, which is also a split quaternion algebra, has been introduced in \cite{cl}.

\section{ Definitions  of the column and  row immanants} Denote by ${\bf H}^{n \times m}$ a set
of $n \times m$ matrices with entries in ${\bf H}$.
For  ${\rm {\bf A}}=(a_{ij})\in{\bf H}^{n\times n}$ we define $n$ row immanents as  follows.
\begin{defn}\label{def_rimm}
The $i$th row immanent of  ${\rm {\bf A}}=(a_{ij})\in{\bf H}^{n\times n}$ is defined by putting
\[
{\rm{rImm}}_{ i} {\rm {\bf A}} = {\sum\limits_{\sigma \in S_{n}}
{\chi(\sigma){a_{i{\kern 1pt} i_{k_{1}}} }
{a_{i_{k_{1}}   i_{k_{1} + 1}}}   \ldots } } {a_{i_{k_{1} + l_{1}}
 i}}  \ldots  {a_{i_{k_{r}}  i_{k_{r} + 1}}}
\ldots  {a_{i_{k_{r} + l_{r}}  i_{k_{r}} }},
\]
where left-ordered cycle notation of the permutation $\sigma$ is written
as follows
\begin{equation}\label{eq:sigma}
   \sigma = \left( {i\,i_{k_{1}}  i_{k_{1} + 1} \ldots i_{k_{1} +
l_{1}} } \right)\left( {i_{k_{2}}  i_{k_{2} + 1} \ldots i_{k_{2} +
l_{2}} } \right)\ldots \left( {i_{k_{r}}  i_{k_{r} + 1} \ldots
i_{k_{r} + l_{r}} } \right).
\end{equation}
Here the index $i$ starts the first cycle from the left  and other
cycles satisfy the following conditions
\begin{equation}\label{eq:order}
 i_{k_{2}}  < i_{k_{3}}  < \ldots < i_{k_{r}},
\quad i_{k_{t}}  < i_{k_{t} + s}.
\end{equation}
for all $t = \overline {2,r}$ and $s = \overline {1,l_{t}}$.
\end{defn}
Consequently we have the following definitions.
\begin{defn}\label{def_rperm}
The $i$th row permanent of ${\rm {\bf A}}=(a_{ij})\in{\bf H}^{n\times n}$ is defined as
\[
{\rm{rper}}_{ i} {\rm {\bf A}} = {\sum\limits_{\sigma \in S_{n}}
{{a_{i{\kern 1pt} i_{k_{1}}} }
{a_{i_{k_{1}}   i_{k_{1} + 1}}}   \ldots } } {a_{i_{k_{1} + l_{1}}
 i}}  \ldots  {a_{i_{k_{r}}  i_{k_{r} + 1}}}
\ldots  {a_{i_{k_{r} + l_{r}}  i_{k_{r}} }},
\]
where left-ordered cycle notation of the permutation $\sigma$ satisfies the conditions (\ref{eq:sigma}) and  (\ref{eq:order}).
\end{defn}
\begin{defn}\label{def_rdet}\cite{ky1}
The $i$th row determinant of  ${\rm {\bf A}}=(a_{ij})\in{\bf H}^{n\times n}$ is defined as
\[
{\rm{rdet}}_{ i} {\rm {\bf A}} = {\sum\limits_{\sigma \in S_{n}}
{\left( { - 1} \right)^{n - r}{a_{i{\kern 1pt} i_{k_{1}}} }
{a_{i_{k_{1}}   i_{k_{1} + 1}}}   \ldots } } {a_{i_{k_{1} + l_{1}}
 i}}  \ldots  {a_{i_{k_{r}}  i_{k_{r} + 1}}}
\ldots  {a_{i_{k_{r} + l_{r}}  i_{k_{r}} }},
\]
where left-ordered cycle notation of the permutation $\sigma$ satisfies the conditions (\ref{eq:sigma}) and  (\ref{eq:order}), (since  ${\rm sign}(\sigma)=\left( { - 1} \right)^{n - r}$).
\end{defn}
For ${\rm {\bf A}}=(a_{ij})\in{\bf H}^{n\times n}$ we  define $n$ column immanents as well.
\begin{defn}\label{def_cimm}
The $j$th column immanent of  ${\rm {\bf A}}=(a_{ij})\in{\bf H}^{n\times n}$ is defined as
\[
{\rm{cImm}}_{ j} {\rm {\bf A}} =  {\sum\limits_{\tau \in S_{n}}
{\chi(\tau)
{a_{j_{k_{r}} j_{k_{r} +
l_{r}} } \ldots a_{j_{k_{r} + 1}  j_{k_{r}} }  \ldots } }a_{j\,
j_{k_{1} + l_{1}} }  \ldots  a_{ j_{k_{1} + 1} j_{k_{1}}
}a_{j_{k_{1}} j}},
\]
where right-ordered cycle notation of the
permutation $\tau \in S_{n}$ is written as follows
\begin{equation}\label{eq:tau}
 \tau = \left( {j_{k_{r} + l_{r}}  \ldots j_{k_{r} + 1} j_{k_{r}}
} \right)\ldots \left( {j_{k_{2} + l_{2}}  \ldots j_{k_{2} + 1}
j_{k_{2}} } \right){\kern 1pt} \left( {j_{k_{1} + l_{1}}  \ldots
j_{k_{1} + 1} j_{k_{1} } j} \right).
\end{equation}
Here the first cycle from the right begins with the index $j$ and
other cycles satisfy the following conditions
\begin{equation}\label{eq:order_tau}
 j_{k_{2}}  < j_{k_{3}}  < \ldots < j_{k_{r}},
\quad j_{k_{t}}  < j_{k_{t} + s},
\end{equation}
for all $t = \overline {2,r}$ and $s = \overline {1,l_{t}}$.
\end{defn}
Consequently we have the following definitions as well.
\begin{defn}\label{def_cperm}
The  $j$th column permanent of  ${\rm {\bf A}}=(a_{ij})\in{\bf H}^{n\times n}$ is defined as
\[
{\rm{rper}}_{ j} {\rm {\bf A}} = {\sum\limits_{\tau \in S_{n}}
{{a_{j_{k_{r}} j_{k_{r} +
l_{r}} } \ldots a_{j_{k_{r} + 1}  j_{k_{r}} }  \ldots } }a_{j\,
j_{k_{1} + l_{1}} }  \ldots  a_{ j_{k_{1} + 1} j_{k_{1}}
}a_{j_{k_{1}} j}},
\]
where right-ordered cycle notation of the permutation $\sigma$ satisfies the conditions (\ref{eq:tau}) and  (\ref{eq:order_tau}).
\end{defn}
\begin{defn}\label{def_cdet}\cite{ky1}
The  $j$th column determinant of  ${\rm {\bf A}}=(a_{ij})\in{\bf H}^{n\times n}$ is defined as
\[
{\rm{rdet}}_{ j} {\rm {\bf A}} = {\sum\limits_{\tau \in S_{n}}
\left( { - 1} \right)^{n - r}{{a_{j_{k_{r}} j_{k_{r} +
l_{r}} } \ldots a_{j_{k_{r} + 1}  j_{k_{r}} }  \ldots } }a_{j\,
j_{k_{1} + l_{1}} }  \ldots  a_{ j_{k_{1} + 1} j_{k_{1}}
}a_{j_{k_{1}} j}},
\]
where right-ordered cycle notation of the permutation $\sigma$ satisfies the conditions (\ref{eq:tau}) and  (\ref{eq:order_tau}).
\end{defn}

\section{Basic properties  of the column and  row immanants}
Consider
the basic properties of the column and row immanents of a
square matrix over ${\rm {\bf{H}}}$.
\begin{prop}\label{theorem:zeros_col_row}(The first theorem about zero of an immanant)
If one of the rows (columns) of  ${\rm {\bf A}}=(a_{ij})\in{\bf H}^{n\times n}$  consists of zeros only, then $
{\rm{rImm}}_{{i}}\, {\rm {\bf A}} = 0$ and ${\rm{cImm}} _{{i}}\,
{\rm {\bf A}} = 0$ for all ${i = \overline {1,n}}. $
\end{prop}
{\textit{Proof}}. The proof immediately
follows from the definitions.
$\blacksquare$

Denote by ${\bf{H}}a$ and $a{\bf{H}}$  left  and  right principal ideals of ${\bf{H}}$, respectively.
\begin{prop}\label{theorem:zeros_row2}(The second theorem about zero of an row immanant)
Let ${\rm {\bf A}}=(a_{ij})\in{\bf H}^{n\times n}$ and  $a_{ki}\in{\bf{H}}a_{i}$ and $a_{ij}\in \overline{a_{i}}{\bf{H}}$, where $n(a_{i})=0$ for $k,j = \overline {1,n}$ and  for all $i\neq k$. Let $a_{11}\in{\bf{H}}a_{1}$ and $a_{22}\in \overline{a_{1}}{\bf{H}}$ if $k=1$, and $a_{kk}\in{\bf{H}}a_{k}$ and $a_{11}\in \overline{a_{k}}{\bf{H}}$ if $k=i>1$, where $n(a_{k})=0$. Then ${\rm{rImm}}_{k} {\rm {\bf A}} =0$.
\end{prop}
{\textit{Proof}}.  Let $i\neq k$. Consider an arbitrary monomial of ${\rm{rImm}}_{k} {\rm {\bf A}}$, if $i\neq k$,
\[
 d = \chi(\sigma)a_{ki }   a_{ij}   \ldots  a_{lm }
\]
where $\{l, m\}\subset \{1,...,n\}$. Since there exists $a_{i}\in {\bf H}$ such that $n(a_{i})=0$, and $a_{ki}\in{\bf{H}}a_{i}$,  $a_{ij}\in \overline{a_{i}}{\bf{H}}$, than $a_{ki } a_{ij}=0$ and $d=0$.

Let $i=k=1$. Then an arbitrary monomial of ${\rm{rImm}}_{1} {\rm {\bf A}}$,
\[
 d = \chi(\sigma)a_{11 }   a_{22}   \ldots  a_{lm }.
\]
Since there exists $a_{1}\in {\bf H}$ such that $n(a_{1})=0$, and $a_{11}\in{\bf{H}}a_{1}$,  $a_{22}\in \overline{a_{i}}{\bf{H}}$, than $a_{11 } a_{22}=0$ and $d=0$.
If $k=i>1$, then an arbitrary monomial of ${\rm{rImm}}_{k} {\rm {\bf A}}$,
\[
 d = \chi(\sigma)a_{kk }   a_{11}   \ldots  a_{lm }.
\]
Since there exists $a_{k}\in {\bf H}$ such that $n(a_{k})=0$, and $a_{kk}\in{\bf{H}}a_{k}$,  $a_{11}\in \overline{a_{k}}{\bf{H}}$, than $a_{kk } a_{11}=0$ and $d=0$.
$\blacksquare$

\begin{prop}\label{theorem:zeros_col2}(The second theorem about zero of an column immanant)
Let ${\rm {\bf A}}=(a_{ij})\in{\bf H}^{n\times n}$ and  $a_{ik}\in a_{i}{\bf{H}}$ and $a_{ji}\in {\bf{H}}\overline{a_{i}}$, where $n(a_{i})=0$ for $k,j = \overline {1,n}$ and  for all $i\neq k$. Let $a_{11}\in a_{1}{\bf{H}}$ and $a_{22}\in {\bf{H}}\overline{a_{1}}$ if $k=1$, and $a_{kk} a_{k}\in{\bf{H}}$ and $a_{11}\in {\bf{H}}\overline{a_{k}}$ if $k=i>1$, where $n(a_{k})=0$. Then ${\rm{cImm}}_{k} {\rm {\bf A}} =0$.
\end{prop}
{\textit{Proof}}. The proof is similar to the proof of Proposition \ref{theorem:zeros_row2}.
$\blacksquare$

The proofs of the next theorems immediately
follow from the definitions.
\begin{prop}\label{theorem:l_mult} If the $i$th row of  ${\rm {\bf A}}=(a_{ij})\in{\bf H}^{n\times n}$ is left-multiplied by $b \in {\bf{H}} $, then $
{\rm{rImm}}_{{i}}\, {\rm {\bf A}}_{{i{\kern 1pt}.}} \left( {b
\cdot {\rm {\bf a}}_{{i{\kern 1pt}.}}} \right) = b \cdot
{\rm{rImm}}_{{i}}\, {\rm {\bf A}}$ for all ${i = \overline {1,n}}
.$
\end{prop}
\begin{prop}\label{theorem:c_mult} If the $j$th column of
 ${\rm {\bf A}}=(a_{ij})\in{\bf H}^{n\times n}$ is right-multiplied by $b \in {\bf{H}}$, then $
{\rm{cImm}} _{{j}}\, {\rm {\bf A}}_{{.{\kern 1pt}j}} \left( {{\rm
{\bf a}}_{{.{\kern 1pt}j}} \cdot b} \right) = {\rm{cImm}} _{{j}}\,
{\rm {\bf A}}\cdot b$ for all ${j = \overline {1,n}}.$
\end{prop}
\begin{prop}\label{theorem:sum_matr_row} If for  ${\rm {\bf A}}=(a_{ij})\in{\bf H}^{n\times n}$ there exists  $t \in \{1,...,n\} $  such that $a_{tj} =
b_{j} + c_{j} $ for all $j = \overline {1,n}$, then for all $i =
\overline {1,n}$
\[
\begin{array}{l}
   {\rm{rImm}}_{{i}}\, {\rm {\bf A}} = {\rm{rImm}}_{{i}}\, {\rm {\bf
A}}_{{t{\kern 1pt}.}} \left( {{\rm {\bf b}}} \right) +
{\rm{rImm}}_{{i}}\, {\rm {\bf A}}_{{t{\kern 1pt}.}} \left( {{\rm
{\bf c}}} \right), \\
  {\rm{cImm}} _{{i}}\, {\rm {\bf A}} = {\rm{cImm}} _{{i}}\, {\rm
{\bf A}}_{{t{\kern 1pt}.}} \left( {{\rm {\bf b}}} \right) +
{\rm{cImm}}_{{i}}\, {\rm {\bf A}}_{{t{\kern 1pt}.}} \left( {{\rm
{\bf c}}} \right),
\end{array}
\]
where ${\rm {\bf b}}=(b_{1},\ldots, b_{n})$, ${\rm {\bf
c}}=(c_{1},\ldots, c_{n}).$
\end{prop}
\begin{prop}\label{theorem:sum_matr_col} If for ${\rm {\bf A}}=(a_{ij})\in{\bf H}^{n\times n}$ there exists  $t \in \{1,...,n\} $ such that
$a_{i\,t} = b_{i} + c_{i}$  for all  $i = \overline {1,n}$, then for all $j
= \overline {1,n}$
\[
\begin{array}{l}
  {\rm{rImm}}_{{j}}\, {\rm {\bf A}} = {\rm{rImm}}_{{j}}\, {\rm {\bf
A}}_{{.\,{\kern 1pt}t}} \left( {{\rm {\bf b}}} \right) +
{\rm{rImm}}_{{j}}\, {\rm {\bf A}}_{{.\,{\kern 1pt} t}} \left(
{{\rm
{\bf c}}} \right),\\
  {\rm{cImm}} _{{j}}\, {\rm {\bf A}} = {\rm{cImm}} _{{j}}\, {\rm
{\bf A}}_{{.\,{\kern 1pt}t}} \left( {{\rm {\bf b}}} \right) +
{\rm{cImm}} _{{j}} {\rm {\bf A}}_{{.\,{\kern 1pt}t}} \left( {{\rm
{\bf c}}} \right),
\end{array}
\]
where ${\rm {\bf b}}=(b_{1},\ldots, b_{n})^T$, ${\rm {\bf
c}}=(c_{1},\ldots, c_{n})^T.$
\end{prop}
\begin{prop}  If ${\rm {\bf A}}^{ *} $
is the Hermitian adjoint matrix (the conjugate transpose) of  ${\rm {\bf A}}=(a_{ij})\in{\bf H}^{n\times n}$, then $
{\rm{rImm}}_{{i}}\, {\rm {\bf A}}^{ *} =  \overline{{{\rm{cImm}}
_{{i}}\, {\rm {\bf A}}}}$ for all $i = \overline {1,n} $.
\end{prop}
Particular cases of these properties for the row-column determinants and permanents  are evident.

\section{ An immanent of a Hermitian matrix}

If ${\rm {\bf A}}^{ *}={\bf A}$ then ${\rm {\bf A}}\in{\bf H}^{n\times n}$ is called a Hermitian matrix. In this section we consider the key theorem about row-column immanats of a Hermitian matrix.

The following lemma is needed for the sequel.
\begin{lem}\cite{ky1} \label{lemma:sum_factors}  Let $T_{n} $ be the sum
of all  possible products  of  the $n$ factors, each of which are
either $h_{i} \in {\bf H}$ or $ \overline {h_{i}} $
for all $ i = \overline {1,n}$, by specifying the ordering in the
terms, $
T_{n} = h_{1} \cdot h_{2} \cdot \ldots \cdot h_{n} + \overline
{h_{1}} \cdot h_{2} \cdot \ldots \cdot h_{n} + \ldots + \overline
{h_{1}}  \cdot \overline {h_{2}}  \cdot \ldots \cdot \overline
{h_{n}}.
$
Then $T_{n}$ consists of the  $2^{n}$ terms and $T_{n} = {\rm
t}\left( {h_{1}} \right)\;{\rm t}\left( {h_{2}} \right)\;\ldots
\;{\rm t}\left( {h_{n}} \right).$
\end{lem}
\begin{thm}\label{theorem:rdet_eq_cdet} If ${\rm {\bf A}}\in{\bf H}^{n\times n}$ is a Hermitian matrix, then
\[
{\rm{rImm}} _{1} {\rm {\bf A}} = \ldots = {\rm{rImm}} _{n} {\rm
{\bf A}} = {\rm{cImm}} _{1} {\rm {\bf A}} = \ldots = {\rm{cImm}}
_{n} {\rm {\bf A}}  \in {\bf{F}}.
\]
\end{thm}
{\textit{Proof}}.  At first we note  that if  ${\rm {\bf A}}=(a_{ij})\in{\bf H}^{n\times n}$ is
Hermitian, then we have $a_{ii} \in{\bf{F}}$ and $a_{ij} =
\overline {a_{ji}}$ for all $i,j = \overline {1,n}$.

 We divide the set of monomials of
${\rm{rImm}}_{i} {\rm {\bf A}}$ for some $i \in \{1,...,n\}$
into two subsets. If indices of coefficients of monomials form
permutations as products of disjoint cycles of length 1 and 2,
then we include these monomials to the first subset. Other
monomials belong to the second subset. If indices of coefficients
form a disjoint cycle of length 1, then these coefficients are
$a_{jj}$ for $j \in \{1,...,n\}$ and $a_{jj}\in {\bf{F}}$.

If indices
of coefficients form a disjoint cycle of length 2, then these
entries are conjugated, $a_{i_{k} i_{k + 1}} = \overline {a_{i_{k
+ 1} i_{k}} } $, and
\[
a_{i_{k} i_{k + 1}}  \cdot a_{i_{k + 1} i_{k}}  = \overline
{a_{i_{k + 1} i_{k}} }  \cdot a_{i_{k + 1} i_{k}}  = {\rm
n}(a_{i_{k + 1} i_{k}}  ) \in {\bf{F}}.
\]
So, all monomials of the first subset take on values in ${\bf{F}}$.

Now we consider some monomial  $d$ of the second subset. Assume
that its index  permutation $\sigma$ forms a direct product of $r$ disjoint
cycles. Denote $i_{k_{1}}:=i$, then
\begin{equation}
\label{kyr6}
\begin{array}{l}
 d = \chi(\sigma)a_{i_{k_{1}}i_{k_{1}+1} }   \ldots  a_{i_{k_{1} + l_{1}
} i_{k_{1}}}  a_{i_{k_{2}}  i_{k_{2} + 1}}   \ldots  a_{i_{k_{2} +
l_{2}} i_{k_{2}} }   \ldots  a_{i_{k_{m}}  i_{k_{m} + 1}} \ldots
\times
\\
 \times a_{i_{k_{m} + l_{m}}  i_{k_{m}} }   \ldots  a_{i_{k_{r}}
i_{k_{r} + 1}}  \ldots  a_{i_{k_{r} + l_{r}}  i_{k_{r}} }  = \chi(\sigma)h_{1}  h_{2}  \ldots  h_{m}  \ldots  h_{r} ,
\\
 \end{array}
\end{equation}
\noindent where $h_{s} = a_{i_{k_{s}} i_{k_{s} + 1}}  \cdot \ldots
\cdot a_{i_{k_{s} + l_{s}}  i_{k_{s}} }$ for all $ s = \overline
{1,r}$, and $ m \in \{1,\ldots ,r\}.$   If $l_{s} = 1$, then
$h_{s} = a_{i_{k_{s}}  i_{k_{s} + 1}} a_{i_{k_{s} + 1{\kern 1pt}}
i_{k_{s}} }  = {\rm n}(a_{i_{k_{s}} i_{k_{s} + 1}} ) \in {\bf{F}}$. If $l_{s} = 0$, then $h_{s} = a_{i_{k_{s}}
i_{k_{s}} } \in {\bf{F}}$. If  $l_{s} =0$ or $l_{s} =1$
for all $s = \overline {1,r}$ in (\ref{kyr6}), then $d$ belongs to the first subset. Let there exists $s \in I_{n} $ such
that $l_{s} \ge 2$. Then
\[
   \overline {h_{s}}  = \overline {a_{i_{k_{s}}  i_{k_{s} + 1}}
\ldots  a_{i_{k_{s} + l_{s}}  i_{k_{s}} } }
   = \overline
{a_{i_{k_{s} + l_{s} } i_{k_{s}} } }   \ldots  \overline
{a_{i_{k_{s}}  i_{k_{s} + 1}} } = a_{i_{k_{s}}  i_{k_{s} + l_{s}}
}  \ldots  a_{i_{k_{s} + 1} i_{k_{s}} }.
\]
Denote by  $\sigma _{s} \left( {i_{k_{s}} }  \right){\rm :} =
\left( {i_{k_{s}}  i_{k_{s} + 1} \ldots i_{k_{s} + l_{s}} }
\right)$
  a disjoint cycle of indices of $d$ for some $ s \in \{1,...,r\}$, then $\sigma=\sigma _{1} \left( {i_{k_{1}} }  \right)\sigma _{2} \left( {i_{k_{2}} }  \right)...\sigma _{r} \left( {i_{k_{r}} }  \right)$.  The disjoint cycle $\sigma _{s} \left(
{i_{k_{s}} }  \right)$ corresponds to the factor $h_{s}$. Then
$\sigma _{s}^{ - 1} \left( {i_{k_{s}} }  \right) = \left(
{i_{k_{s}} i_{k_{s} + l_{s}} i_{k_{s} + 1} \ldots i_{k_{s} + 1}}
\right)$ is the inverse disjoint cycle and $\sigma _{s}^{ - 1}
\left( {i_{k_{s}} } \right) $ corresponds to the factor
$\overline{h_{s}}$. By Lemma \ref{lemma:sum_factors} there exist
another $2^{p}-1$ monomials for $d$, (where $p = r - \rho $ and
$\rho $ is the number of disjoint cycles of length 1 and 2), such
that their index permutations form the direct products of $r$
disjoint cycles either $\sigma _{s} \left( {i_{k_{s}} } \right)$
or $\sigma _{s}^{ - 1} \left( {i_{k_{s}} } \right)$ by specifying
their ordering by $s$ from $1$ to $r$. Their cycle notations are
left-ordered according to Definition \ref{def_rimm}. These permutations are unique decomposition of the permutation $\sigma$ including their ordering by $s$ from $1$ to $r$.  Suppose $C_{1} $
is the sum of these $2^{p}-1$ monomials and $d$, then by Lemma
\ref{lemma:sum_factors} we obtain
\[
C_{1} = \chi(\sigma)\alpha \;{\rm t}(h_{\nu _{1}}  )\;\ldots
\;{\rm t}(h_{\nu _{p}}  ) \in {\bf{F}}.
\]
\noindent Here $\alpha \in {\bf{F}}$ is the product of
coefficients
 whose indices   form disjoint cycles of length 1 and 2,
  $\nu _{k} \in \{1,\ldots ,r\}$ for all $ k = \overline {1,p}$.

Thus for an arbitrary monomial of the second subset of
${\rm{rImm}}_{i}\, {\rm {\bf A}}$, we can  find the $2^{p}$
monomials  such that their sum  takes on a value in ${\bf{F}}$. Therefore, ${\rm{rImm}} _{i}\, {\rm {\bf A}}\in
{\bf{F}}$.

Now we  prove the equality of all  row immanents of ${\rm {\bf
A}}$. Consider an arbitrary  ${\rm{rImm}} _{j}\, {\rm {\bf A}}$
such that $j \ne i$ for all ${\ j = \overline {1,n}}$. We divide
the set of monomials of ${\rm{rImm}} _{j}\, {\rm {\bf A}}$ into
two subsets using the same rule as for ${\rm{rImm}}_{i}\, {\rm
{\bf A}}.$ Monomials of the first subset are products of entries
of the principal diagonal or  norms of entries of ${\rm {\bf A}}$.
Therefore they  take on a value in ${\bf{F}}$ and each
monomial of the first subset of ${\rm{rImm}}_{i}\, {\rm {\bf A}}$
is equal to a corresponding monomial of the first subset of
${\rm{rImm}} _{j}\, {\rm {\bf A}}$.

Now consider the monomial $d_{1} $ of the second subset of
monomials of ${\rm{rImm}}_{j}\, {\rm {\bf A}}$ consisting of
coefficients that are equal to the coefficients of $d$ but they are in another order.
 Consider all  possibilities
of the arrangement of coefficients in $d_{1} $.

 (i) Suppose that the index permutation ${\sigma}'$ of its coefficients form
 a  direct product of $r$ disjoint cycles and these cycles coincide
  with the $r$ disjoint cycles of $d$ but differ by their ordering.
 Then ${\sigma}'=\sigma$ and we have
\[d_{1} = \chi(\sigma)\alpha h_{\mu}   \ldots  h_{\lambda}  ,
\]
where $\{\mu ,\ldots ,\lambda \} = \{\nu _{1} ,\ldots ,\nu _{p}
\}$. By Lemma \ref{lemma:sum_factors} there exist $2^{p }- 1$
monomials  of the second subset of ${\rm{rImm}} _{j}\, {\rm {\bf
A}}$ such that each of them is equal to a product of $p$ factors
either $h_{s} $ or $\overline {h_{s}} $ for all $s \in \{\mu
,\ldots ,\lambda \}$. Hence
by Lemma \ref{lemma:sum_factors}, we obtain
\[
C_{2} = \chi(\sigma)\alpha \;t(h_{\mu}  )\;\ldots
\;t(h_{\lambda}  ) = \chi(\sigma)\;\alpha \;t(h_{\nu _{1}}
)\ldots \;t(h_{\nu _{p}}  ) = C_{1} .
\]
 (ii) Now suppose that  in addition to the case (i)  the index $j$
 is
placed inside some disjoint cycle of the index permutation $\sigma$ of $d$,
e.g. $ j \in \{i_{k_{m}+1},...,i_{k_{m}+l_{m}}\}$. Denote $j =
i_{k_{m} + q} $. Considering the above said and $\sigma_{k_{m}+1} (i_{k_{m}+1})=\sigma_{k_{m}+q} (i_{k_{m}+q})$, we have ${\sigma}'=\sigma$. Then $d_{1} $ is represented as follows:
\begin{equation}\label{kyr7}
\begin{array}{c}
   d_{1} = \chi(\sigma)a_{i_{k_{m} + q} i_{k_{m} + q + 1}}  \ldots
\quad  a_{i_{k_{m} + l_{m}}  i_{k_{m}} }  \, a_{i_{k_{m}} i_{k_{m}
+ 1}} \ldots   \times \\
   \times   a_{i_{k_{m} + q - 1} i_{k_{m} + q}} a_{i_{k_{\mu} }  i_{k_{\mu}  + 1}}
      \ldots  a_{i_{k_{\mu
} + l_{\mu} }  i_{k_{\mu} } }   \ldots a_{i_{k_{\lambda} }
i_{k_{\lambda}  + 1}}  \ldots a_{i_{k_{\lambda}  + l_{\lambda} }
i_{k_{\lambda} } }  = \\
  =\chi(\sigma)\alpha  \tilde {h}_{m}  h_{\mu}
\ldots  h_{\lambda},
\end{array}
\end{equation}
\noindent where $\{m,\mu ,\ldots ,\lambda \} = \{\nu _{1} ,\ldots
,\nu _{p} \}$. Except for $\tilde {h}_{m} $,  each factor of
$d_{1}$ in (\ref{kyr7}) corresponds to the equal factor of $d$ in
(\ref{kyr6}). By the rearrangement property of the trace, we have
$t(\tilde {h}_{m} ) = t(h_{m} )$. Hence by Lemma
\ref{lemma:sum_factors} and by analogy to the previous case, we
obtain,
\[
\begin{array}{c}
  C_{2} = \chi(\sigma)\alpha \;t(\tilde {h}_{m} )\;t(h_{\mu}
)\;\ldots \;t(h_{\lambda}  )=\\
 = \chi(\sigma)\;\alpha \;t(h_{\nu _{1}}  )\;\ldots \;t(h_{m} )\;\ldots
\;t(h_{\nu _{p}}  ) = C_{1} .
\end{array}
\]
(iii) If in addition to the case (i)  the index $i$ is placed
inside some disjoint cycles of the index permutation of $d_{1}$,
then we apply the rearrangement property of the trace  to this
cycle. As in the previous cases  we find  $2^{p}$ monomials of the
second subset of ${\rm{rImm}} _{j} \,{\rm {\bf A}}$ such that by
Lemma \ref{lemma:sum_factors} their sum is equal to the sum of the
corresponding $2^{p}$ monomials of ${\rm{rImm}}_{i} {\rm {\bf
A}}$. Clearly, we obtain the same conclusion at association of all
previous cases, then we apply twice the rearrangement property of
the trace.

Thus, in any case  each  sum of  $2^{p}$ corresponding monomials
of the second subset of ${\rm{rImm}} _{{j}}\, {\rm {\bf A}}$ is
equal to the sum of  $2^{p}$ monomials of ${\rm{rImm}}_{{i}}\,
{\rm {\bf A}}$. Here $p$ is the number of disjoint cycles of
length  more than 2. Therefore, for all $i, j = \overline {1,n}$
we have
\[
{\rm{rImm}}_{{i}}\, {\rm {\bf A}} = {\rm{rImm}} _{{j}}\, {\rm {\bf
A}}  \in {\bf{F}}.
\]
The equality ${\rm{cImm}} _{{i}}\, {\rm {\bf A}} =
{\rm{rImm}}_{{i}}\, {\rm {\bf A}}$ for all $ i = \overline {1,n}$ is proved similarly.
$\blacksquare$

Since Theorem \ref{theorem:rdet_eq_cdet} we have the following definition.
\begin{defn}\label{def_det__her}
Since all  column and  row immanents of a
Hermitian matrix over ${\bf{H}}$ are equal, we can
define the immanant (permanent, determinant) of a  Hermitian matrix ${\rm {\bf A}}\in{\bf H}^{n\times n}$. By definition, we
put for all $i = \overline {1,n}$
\[\begin{array}{c}
  {\rm Imm}\, {\rm {\bf A}}: = {\rm{rImm}}_{{i}}\,
{\rm {\bf A}} = {\rm{cImm}} _{{i}}\, {\rm {\bf A}},\\
   {\rm per} \,{\rm {\bf A}}: = {\rm{rper}}_{{i}}\,
{\rm {\bf A}} = {\rm{cper}} _{{i}}\, {\rm {\bf A}},\\
   \det {\rm {\bf A}}: = {\rm{rdet}}_{{i}}\,
{\rm {\bf A}} = {\rm{cdet}} _{{i}}\, {\rm {\bf A}}.
\end{array}
\]
\end{defn}

\end{document}